\documentclass[11pt]{amsart}

\newcommand{\R}{\mathbb R}%
\newcommand{\I}{\mathbb I}%
\newtheorem{theorem}{Theorem}[section]

\begin{document}

\title[Gaussian noise measure]{Nash twist and Gaussian noise measure for\\ isometric $C^1$ maps}
\author[A. Dasgupta]{Amites Dasgupta}
\address{Statistics and Mathematics Unit, Indian Statistical Institute\\ 203,
B.T. Road, Calcutta 700108, India.\\ e-mail:amites@isical.ac.in\\ }
\author[M. Datta]{Mahuya Datta}
\address{Statistics and Mathematics Unit, Indian Statistical Institute\\ 203,
B.T. Road, Calcutta 700108, India.\\ e-mail:mahuya@isical.ac.in\\ }
\begin{abstract}Starting with a short map $f_0:I\to \mathbb R^3$ on the unit interval $I$, we
construct random isometric map $f_n:I\to \mathbb R^3$ (with respect to some
fixed Riemannian metrics) for each positive integer $n$, such that the
difference $(f_n - f_0)$ goes to zero in the $C^0$ norm. The construction of
$f_n$ uses the Nash twist. We show that the distribution of $ n^{3/2} (f_n -f_0)$ converges (weakly) to a Gaussian noise measure.
\end{abstract}
\maketitle
\section{introduction}

The problem of associating a measure to the solution space of a differential equation has been mentioned by Gromov in an interview with M. Berger \cite{berger}.
Our point of interest lies in the space of isometric immersions of a Riemannian manifold $(M,g)$ into a Euclidean space $\R^q$ with the canonical metric $h$.
In 1954, Nash proved that if a manifold $M$ with a Riemannian metric $g$ can be embdedded in a Euclidean space $\R^q$, $q > n+1$, then one can construct a large class of isometric $C^1$ embeddings (\cite{nash}). If the initial embedding $f_0:M\to \R^q$ is strictly short, that is if $g-f^*h$ is a Riemannian metric on $M$ then the isometric embeddings can be made to lie  in an arbitrary $C^0$ neighbourhood of the initial embedding. In the following year, Kuiper showed that the bound can be improved to $q\geq n+1$ (\cite{kuiper}). The Nash process is an iterative process;
each stage of the iteration consists of several small steps each of which involves a choice of a rapidly oscillating function defining a perturbation, called a Nash twist.
Starting with the short map $f_0$, one constructs a sequence $\{f_n\}$ of short immersions, where $f_n$ is obtained from $f_{n-1}$ possibly through infinite steps each involving a Nash twist. Successive Nash twists performed on $f_n$ results in a correction to the induced metric $f_n^*h$.
These corrections do not yield an isometric immersion at any stage but each $f_n$ still remains strictly short; however, $f_{n+1}$ is better than $f_n$ in the sense that the induced metric $f_{n+1}^*h$ is closer to $g$ than that in the previous stage. The Nash twist is a controlled perturbation - the $C^1$ distance between any two consecutive maps $f_n$ and $f_{n+1}$ remains bounded by the distance between $g$ and the induced metric $f_n^*h$; furthermore $f_n$ can be made to lie in an arbitrary $C^0$ neighbourhood of $f_0$. As a result the sequence converges to a $g$-isometric $C^1$ immersion.

Nash-Kuiper theory was later generalised by Gromov into the theory of convex integration (\cite{gromov}). Camillo De Lellis and L\'{a}szl\'{o} Sz\'{e}kelyhidi, Jr., have briefly mentioned about the probabilistic approach to convex integration (\cite{lellis}) by pointing to the fact that Convex integration can be seen as a control problem: at each step of the iteration,
one has to choose an admissible perturbation, consisting essentially of a (plane-)wave direction and a frequency.

In the present article we shall consider the domain space of the maps to be 1-dimensional.
In dimension 1, the solution to the $C^1$-isometric embedding problem does not require an infinite Nash process. Isometric maps  $f:\mathbb I\to \R^3$ with respect to the standard Riemannian metrics can be obtained simply by integrating a curve in the 2-sphere and this reduces the Nash process to a single stage. However, Nash twists play an important role in controlling the distance between the initial and the perturbed map - it can produce isometric immersions within an arbitrary $C^0$ neighbourhood of the initial embedding $f_0$. In order to keep the solutions sufficiently close to the original embedding, Nash introduced a periodic function of high frequency (or rapidly oscillating function) under the integration process. We may remark here that in higher dimension, each step in the Nash process can be reduced to a parametric version of the 1-dimensional Nash process described above. The $C^0$-closeness will then translate into $C^\perp$-closeness (refer to \cite[Pp. 170]{gromov}). However,
the problem in dimension greater than 1 is considerably more difficult and we plan to take it up in future.

It is indeed the case that as the frequency in the Nash twist goes to infinity the distance between the initial short map and the resulting isometric immersion goes to zero in the $C^0$-norm. This motivates us to study
the distribution of $f-f_0$ with respect to an appropriate measure on the space of isometric immersions $f:\I\to \R^3$.

We naturally incorporate a randomness in the Nash twist which translates into a Gaussian noise measure for the difference function $f-f_0$. For each positive integer $n$, we construct random functions $f_n$ such that the difference $(f_n(.) - f_0(.))$ goes to zero (in $C^0$ norm). We scale
it up and examine the distribution. We show that the distribution of  $ n^{3/2} (f_n - f_0)$ converges (weakly) to a Gaussian noise measure. Thus the random
solutions $f_n(.)$ can be thought of as distributed like $f_0+n^{-3/2}$(Gaussian noise) for large $n$. In Theorem~\ref{main} we state this rigorously
identifying the weak limit of $n^{3/2}\int_0^t(f_n - f_0)(s)\,ds$ as a Gaussian process. Section 2.2 is devoted to the proof which requires essentially weak convergence of random walks. In the last section, we compare the above process with a class of extensively studied Gaussian processes.

\section{Notation and main result}
Let $M$ be a smooth manifold with a Riemannian metric $g$. The isometric immersions $f:M\to \R^q$ are solutions to the following system of partial differential equations:
\[\langle\frac{\partial f}{\partial u_i},\frac{\partial f}{\partial u_j}\rangle = g_{ij},\ \ \ i,j=1,2,\dots,n,\]
where $u_1,u_2,\dots u_n$ is a local coordinate system on $M$, $g_{ij}$, $i,j=1,2,\dots,n$, are the matrix coefficients of $g$ and the $\langle\ ,\ \rangle$ denotes the inner product on $\mathbb R^q$.

To motivate the concept of randomness and measure we consider the following simple case where the domain space $M$ is the unit interval $[0,1]$ and hence an arbitrary metric on $M$ is of the form $g\,dt^2$, where $g : [0, 1] \rightarrow \mathbb{R}_+$ is a smooth positive function on $[0,1]$. The shortness condition on a smooth regular curve $f_0 : [0, 1] \rightarrow
\mathbb{R}^3$ then translates into the pointwise inequality $0<||\partial_u f_0|| < \sqrt{g}$ on $[0,1]$. Given such an $f_0$ we want to find a function $f_n
: [0, 1] \rightarrow \mathbb{R}^3$ such that $||\partial_u f_n|| = \sqrt{g}$ (which means that $f_n$ is isometric)
and the $C^0$-distance between $f_n$ and $f_0$ decreases with $n$.

The following considerations illustrate the \emph{Nash twist} in dimension 1.
Suppose that $(X,Y,Z)$ is the Frenet-Serret frame along $f_0$ (assuming that such a frame exists at all points $u\in [0,1]$), where $X$ is
the unit tangent along the curve. Then $(Y,Z)$ span a plane field $J$ along $f_0$ perpendicular to $X$. Consider the curve
$Y(u)\cos 2\pi s+Z(u)\sin 2\pi s$, $0 \leq s \leq 1$, on $J(u)$ for each fixed $u\in [0,1]$. Then
with $r^2 = g - ||\partial_u f_0||^2$ the function
$$ \partial_u f_0 + r(u) (Y(u)\cos 2\pi s+Z(u)\sin 2\pi s)$$ has the required euclidean norm $\sqrt{g}$ and over $s \in [0, 1]$
integrates to $\partial_u f_0$ (the convex integration condition). Now, a Nash twist of $f_0$ is given by
$$f_n(t) = f_0(0) + \int_0^t
\{ \partial_u f_0(u) +  r(u) (Y(u)\cos 2\pi nu+Z(u)\sin 2\pi nu) \} du,$$
where $n$ connects with the frequency of the periodic functions, namely $\cos 2\pi nu$ and $\sin 2\pi nu$, mentioned in the previous section.
Clearly, $f_n$ is a solution of the isometry equation since $(Y(u),Z(u))$ is an orthonormal basis of $J(u)$.
We want to show that the function (or the difference curve) $\int_0^t r(u) (Y(u)\cos 2\pi nu+Z(u)\sin 2\pi nu) du$ is uniformly small over $[0, 1]$. We do this for the two integrals separately with some notational abuse.

Applying integration by parts we get
\begin{eqnarray*} \int_0^1 r(u) \cos {2\pi nu}\, du
&=& \frac{1}{2\pi n}r(u)\sin 2\pi nt - \frac{1}{2\pi n}\int_0^t r'(u)\sin 2\pi nu\, du
\end{eqnarray*}
which in absolute value is bounded by const.$\frac{1}{2n\pi}$ as $r$ is a smooth function on the interval $[0,1]$.
The same estimates also apply to $\int_0^t r(u)\sin {2\pi nu}\, du$ and thus we conclude uniform closeness of $f_n$ and $f_0$.

This difference curve can also be considered for a random path
by changing the function $H_n(u) = nu$ over random choices.
The $H_n$ in the above example can be obtained by integrating the constant function $h_n=n$ over $[0,1]$.
Instead we take $h_n=\pm n$ on each subinterval $(k/n, (k+1)/n]$ independently
with equal probability.
These choices are actually explicitly mentioned
by Gromov except for the probability part. Then integration of the resulting function
will give a random function $H_n(\omega, u)$ which is the graph of a simple random
walk. To see this calculation we consider a sequence of independent and identically distributed random variables $X_k$ which take the values $\pm 1$ with equal
probability on a probability space $(\Omega,\mathcal F, P)$, where $\Omega$ can be taken as the infinite product space $\{-1,+1\}^{\mathbb N}$ and
consider the function $h_n(\omega, x) = nX_k, (k-1)/n \leq x < k/n$. Therefore, each subinterval of length $1/n$ contributes $\pm 1$ and hence
\[H_n(\omega, u) =
\int_0^u h_n(\omega, x) dx = S_k \pm n(u - (k/n)), \ \ k/n \leq u < (k+1)/n,\] where $S_k = X_1 +
\cdots + X_k$. The random sum $S_n(\omega, t)$ can be interpreted as the random walk obtained by linearly joining $S_k$, $1\leq k\leq n$.

If we now consider the components of the random difference curve
\[\int_0^t r(u) [Y(u) \cos {2\pi H_n(\omega, u)}+Z(u) \sin {2\pi H_n(\omega, u)}] du,\]
$C^0$-closeness will follow in the same way. However from the viewpoint of probability,
when a sequence of random variables $V_n$ converges almost surely to a random variable
$V$, many times the difference $V_n - V$, after scaling, converges in a suitable sense
to a nontrivial limiting random variable. When the convergence is weak convergence,
the resulting distribution of the limit is a measure, on a suitable space, associated to the sequence (\cite{billingsley}).

In our case the graph
of $\int_0^t e^{2\pi iH_n(\omega, u)} du$ (omitting $r(u), Y(u), Z(u)$ for simplicity now)
looks similar (in a probabilistic sense
of considering all possible paths) over equal intervals and
is independent over disjoint intervals. However, the limit of
its normalization is not a function but a random distribution (in the sense of generalized functions),
it is called the Brownian white noise measure. To keep track of the random difference curve one
tracks its rescaled integral, which converges weakly to Brownian motion. A rigorous
formulation (bringing in $r(u)$, $Y(u)$, $Z(u)$) is the following

\begin{theorem} The sequence of processes
$2 \pi n^{3/2} \int_0^t (f_n-f_0) \,ds, 0 \leq t \leq 1$,
converges weakly to
$$\int_0^t r(u) Z(u) dW(u) - \int_0^t \int_0^s \partial_u (r(u) Z(u)) dW(u)\, ds, 0 \leq t \leq 1,$$ as
$n \rightarrow \infty$, where $W(\cdot)$ denotes the Wiener measure on $C[0, 1]$, the space of real valued continuous functions on the interval $[0,1]$.
\label{main}\end{theorem}
In this sense, the limit of the scaled random difference curve is locally the Gaussian noise
measure $r(t) Z(t) dW(t) - (\int_0^t \partial_u (r(u) Z(u)) dW(u))\, dt$ (see \cite{mitoma}).
The rate of convergence may also be interesting from the probabilistic point of view.

\section{Proof of the main result}

We first keep track of the integral (omitting $r(u), Y(u), Z(u)$)
\begin{equation}\label{wn}
\int_0^s e^{2 \pi i H_n(u)} du = \int_0^s e^{2 \pi i n u } du =
(1/2\pi n)[  \sin (2 \pi ns)  - i \{\cos (2 \pi n s) -
1 \}], 0 \leq s \leq 1.
\end{equation} The main observation
about this function on the right  is that over each $(k/n, (k+1)/n]$ interval
the imaginary part is the graph of $1 - \cos (2 \pi n x), 0 < x \leq 1/n$, and
the real part is the graph of $\sin (2 \pi n x), 0 < x \leq 1/n$. The
periodic behavior along with the factor $(1/2\pi n)$ indicates the
$C^0$-closeness as $s$ varies in $[0, 1]$.

Plugging in $H_n(\omega, u)$ described in the previous section
the corresponding integral $\int_0^s e^{2 \pi i H_n(\omega, u)} du$  differs from (\ref{wn})
in randomly inverting the imaginary part of the graph over each $(k/n, (k+1)/n]$ interval. To see this
consider $s \in (k/n, (k+1)/n]$. Then
$$ \int_{\frac{k}{n}}^s e^{2 \pi i (S_k/n \pm n(u - k/n)} du
= (1/2\pi n)[ \sin (2 \pi n(s - k/n))  \pm (-i) \{\cos (2 \pi n (s - k/n)) - 1 \}  ]. $$
Noting that the integral over each $(k/n, (k+1)/n]$
is zero and using the periodicity of $\sin$ and $\cos$ functions we have established
that the graph of $\int_0^s e^{2 \pi i H_n(\omega, u)} du$ is obtained by randomly inverting the imaginary
part of the graph of (\ref{wn}) over each $(k/n, (k+1)/n]$ interval.
We now prove  that
\[2 \pi n^{3/2} \int_0^t \{ \int_0^s \sin{2 \pi H_n(\omega, u)} du\} ds,\ 0 \leq t \leq 1,\]
converges weakly to Brownian motion as $n \rightarrow \infty$
and
\[2 \pi n^{3/2} \int_0^t \{ \int_0^s \cos{2 \pi H_n(\omega, u)} du\} ds,\ 0 \leq t \leq 1,\]
converges to the zero process.

To see the exact form of this (random) function we note that, as proved, the
function $\int_0^s e^{2 \pi i H_n(\omega, u)} du$ has a graph which
is the graph of (\ref{wn}) with the imaginary part randomly inverted over each $[k/n, (k+1)/n)$
interval. The integral of the (periodic) function in (\ref{wn}) over each $[k/n, (k+1)/n)$
interval is $ (i/2\pi n^2)$ (the real part integrates to contribute zero).
Thus the integral of the function $\int_0^s e^{2 \pi i H_n(\omega, u)} du$
is $\pm i (1/2\pi n^2)$ over the same interval, the $\pm$ sign coming from the
random inverting. Written explicitly, for $k/n \leq t < (k+1)/n$,
\[\begin{array}{rcl}\int_0^t\{\int_0^s e^{2 \pi i H_n(\omega, u)} du\}ds & = &  \sum_{j = 0}^{k - 1}\int _{j/n}^{(j+1)/n}\{\int_0^s e^{2 \pi i
H_n(\omega, u)} du\}ds\\ && + \int _{k/n}^{t}\{\int_0^s e^{2 \pi i
H_n(\omega, u)} du\}ds \\
& = &  \frac{i}{2\pi n^2}\sum_{j = 1}^{k} X_j + O(\frac{1}{n^2}),\end{array}\]
where $X_j$ are independent $\pm 1$ random variables.
Note that the integral from $k/n$ to $t$ adds a continuous function of order $O(\frac{1}{n^2})$ to the random walk obtained
from the $X_j$'s. Multiplying by $2 \pi n^{3/2}$ we get the weak convergence to Brownian motion. In this sense
the random difference $\int_0^s e^{2 \pi i H_n(\omega, u)} du, 0 \leq s \leq 1,$
when normalized converges to the generalized derivative of Brownian motion,
called Brownian white noise.

For the general case, we
consider (with some abuse of notation)
\[\int_0^s r(u) e^{2 \pi i H_n(\omega, u)} du, 0 \leq s \leq 1,\]
where $r$ is a sufficiently differentiable function. In this case, depending on $r$, after two integrations  we get back a different
random walk minus the area under another random walk, and consider weak convergence again.

We shall deal with the real and the imaginary part of the integral separately. As observed before the randomness has no role to play in the real part.
A straightforward calculation shows that

\begin{center}
$\begin{array}{rcl}\int_0^s r(u) \cos(2 \pi H_n(\omega, u))du & = & \int_0^s r(u) \cos(2 \pi nu)du\\
& = &\frac{1}{2n\pi}r(s)\sin 2n\pi s - \frac{1}{2n\pi}\int_0^s r'(u)\sin 2n\pi u \,du\\
& = &\frac{1}{2n\pi}r(s)\sin 2n\pi s + \frac{1}{4n^2\pi^2}r'(s)\cos 2n\pi s\\
& - & \frac{1}{4n^2\pi^2}r'(0)- \frac{1}{4n^2\pi^2}\int_0^s r''(u)\cos 2n\pi u \,du
\end{array}$\end{center}
Therefore,
\begin{center}$\int_0^t\{\int_0^s r(u) \cos(2 \pi H_n(\omega, u))du\}\,ds = \frac{1}{2\pi n}\int_0^t r(s) \sin (2 \pi ns)\, ds + O(\frac{1}{n^2})$.\end{center}
Since the first term on the right hand side is $O(1/n^2)$ it follows that
\[\lim_{n\to\infty} n^{3/2}\int_0^t\{\int_0^s r(u) \cos{2 \pi H_n(\omega, u)} du\}\,ds=0.\]
Thus the real part of the integral when scaled by $n^{3/2}$ converges to zero uniformly.

To deal with the imaginary part of the integral, for $l/n\leq s<(l+1)/n$, we split it as follows
\begin{eqnarray}\label{area}
\sum_{k = 0}^{[ns] - 1} \int_{k/n}^{(k+1)/n} r(u)  \sin (2\pi H_n(\omega, u)) du + \int_{l/n}^s r(u)  \sin (2\pi H_n(\omega, u)) du,
\end{eqnarray}
and then writing $r(u) = r(k/n) + (r(u) - r(k/n))$ on the subinterval $(k/n, (k+1)/n]$  we get the following for (\ref{area}):

\begin{eqnarray}\label{area2}
&  & \sum_{k = 0}^{[ns] - 1} \int_{k/n}^{(k+1)/n} (r(u) - r(k/n))  \sin (2\pi H_n(\omega, u)) du \nonumber \\
&+&r(l/n)\int_{l/n}^s  \sin(2\pi H_n(\omega,u)) du + \int_{l/n}^s (r(u) - r(l/n))  \sin (2\pi H_n(\omega, u)) du.
\end{eqnarray}
We denote the function (represented by the sum) on the first row by $\psi_1(s)$ and the two terms on the second row by $\psi_2(s)$ and $\psi_3(s)$ respectively. Since $n^{3/2}\int_0^t \psi_2(s) ds$ gives a random walk with steps $\pm\frac{1}{\sqrt{n}}r(k/n)$ over the interval $[k/n,(k+1)/n)$, $n^{3/2}\int_0^t \psi_2(s) ds$ converges weakly to $\frac{1}{2\pi}\int_0^t r(s) dW(s)$.

Next we consider the part $\psi_1$. The summands of $\psi_1$ are not necessarily zero as $r$ is non-constant. Substituting $u=k/n+z/n$ in the $k$-th summand and disregarding the $\pm$ signs,
we get $$\frac{1}{n} \int_0^1 [r\left(\frac{k}{n} + \frac{z}{n}\right) - r\left(\frac{k}{n}\right)] \sin (2 \pi z) dz. $$
To clearly understand this contribution
divide each integral $[k/n,(k+1)/n)$ into the parts where the sine function has the same sign (to apply the mean value theorem
for integrals for some  $z_1 \in (0, 1/2)$). Thus,
\begin{eqnarray*}
&& \frac{1}{n} \int_0^1 [r\left(\frac{k}{n} + \frac{z}{n}\right) - r\left(\frac{k}{n}\right)] \sin (2 \pi z) dz \\
&=&  \frac{1}{n}  \int_0^{1/2} [r\left(\frac{k}{n} + \frac{z}{n}\right) - r\left(\frac{k}{n}\right)] \sin (2 \pi z) dz \\
& + &  \frac{1}{n}\int_{1/2}^1 [r\left(\frac{k}{n} + \frac{z}{n}\right) - r\left(\frac{k}{n}\right)] \sin (2 \pi z) dz \\
&=& \frac{1}{n} \int_0^{1/2} [r\left(\frac{k}{n} + \frac{z}{n}\right) - r\left(\frac{k}{n} + \frac{z}{n}+\frac{1}{2n}\right)] \sin (2 \pi z)\, dz \\
&=& \frac{1}{n\pi} \left[r\left(\frac{k}{n} + \frac{z_1}{n}\right) - r\left(\frac{k}{n} + \frac{z_1}{n}+\frac{1}{2n}\right)\right] \\
&=& - \frac{1}{2n^2\pi} r^{\prime} \left(\frac{k}{n} + \frac{z_1}{n}+\frac{\theta}{2n}\right),
\end{eqnarray*}
where $0<\theta<1$. If we add these integrals after multiplying each them by $\pm 1$ from random inversions, and scale the sum by $n^{3/2}$ then it corresponds to a random walk converging weakly to
$$ - \frac{1}{2\pi} \int_0^s r^\prime(u) dW(u).$$
The continuous random curve $\psi_1+\psi_3$ matches this random walk at the points $k/n$ and is otherwise at a distance at most
$$O(\frac{1}{n}|r(k/n + 1/2n) - r(k/n)|) = O(1/n^2)$$ from it. Thus after scaling by $n^{3/2}$ the random curve $\psi_1+\psi_3$  converges
weakly to the same limit and the integral of $\psi_1+\psi_3$
converges weakly to $-\frac{1}{2\pi}\int_0^t \{\int_0^s r^\prime(u) dW(u)\} ds$ by the continuous mapping theorem. This completes the proof of the main result.
For completeness we indicate how the random walk
\begin{equation}\sum_{k=1}^{[nt_1]} r(k/n)Y_k\label{int_psi_2}\end{equation}
(here $Y_k$ are i.i.d. $\pm\frac{1}{\sqrt{n}}$ random variables) and the area under the random walk
\begin{equation}- \sum_{i=1}^{[ns]} r'(i/n)Y_i\label{area_random}\end{equation}
up to time $t_2$ converges jointly in distribution, after which tightness on product space can be used to conclude weak convergence on $C[0,1]\times C[0,1]$. To obtain the area under (\ref{area_random}), each random walk height is multiplied by $1/n$ and then added. The interchange of summation gives the area as
\begin{equation}-\sum_{i=0}^{k-1} (\frac{k-i}{n}) r'(i/n) Y_i,\label{int_psi_1+psi_3}\end{equation}
where $k=[nt_2]$. From (\ref{int_psi_2}) and (\ref{int_psi_1+psi_3}) the limiting joint finite dimensional distribution follows. The limit of the expression in (\ref{int_psi_1+psi_3}) is seen to be $-\frac{1}{2\pi}\int_0^t(t-u)r'(u)\, dW(u)$ which equals $-\frac{1}{2\pi}\int_0^t \{\int_0^s r^\prime(u) dW(u)\} ds$. Now the sum of the two processes (\ref{int_psi_2}) and (\ref{int_psi_1+psi_3}) converges weakly.
\hfill $\Box$

\section{Concluding remarks:}
It is seen from the proof of Theorem~\ref{main} that the random Nash twist on the initial short curve $f_0$ to obtain an increase of $r^2(u)\,du^2$ to the induced metric
$f_0^*h$, leads to a process whose structure is similar to the
following process (refer to \cite[Theorem 6.3]{hida})
$$X(t)=W(t)-\int_0^t\int_0^s\ell(s,u)dW(u)\,ds,$$
where $\ell(s,u)$ is a Volterra kernel with appropriate conditions.
In our case, componentwise we need a Gaussian process $\int_0^t r(u)\,dW(u)$ and a function $r'(u)/r(u)$ to replace $W(t)$ and  $\ell(s,u)$ in the above formula.

In higher dimension the difference metric $g-f_0^*h$ can be written as a finite sum of monomials $r^2 d\varphi^2$, where $\varphi$ is a rank 1 function. Applying Nash twist along (Im\,$Df_0)^\perp$, one is able to add $r^2 d\varphi^2$ only approximately.
To look into the problem of $C^0$ distance one may need several independent Brownian motions in different directions and we refer to Theorem 3.2.5 of Kallianpur and Xiong (\cite{kallianpur}) for such a construction. However a precise formulation combining various directions is not clear and we hope to explore these aspects in future.

\emph{Acknowledgement}: The second author is greatly indebted to Misha Gromov for sharing his insight on the subject during a visit of the author to IHES.

\end{document}